\documentclass[11pt,a4paper]{article}
\usepackage[T1]{fontenc}
\usepackage[utf8]{inputenc}
\usepackage[english]{babel}
\usepackage{lmodern,microtype,amsmath,amssymb,amsthm,mathtools}
\usepackage[margin=26mm]{geometry}
\usepackage{enumitem,tikz-cd,xurl,needspace}
\usepackage[colorlinks=true,linkcolor=blue!45!black,citecolor=blue!45!black,urlcolor=blue!45!black]{hyperref}
\setlist{itemsep=2pt,topsep=5pt}
\allowdisplaybreaks[2]
\newtheorem{theorem}{Theorem}[section]
\newtheorem{proposition}[theorem]{Proposition}
\newtheorem{lemma}[theorem]{Lemma}
\newtheorem{corollary}[theorem]{Corollary}
\theoremstyle{definition}
\newtheorem{definition}[theorem]{Definition}
\newtheorem{example}[theorem]{Example}
\theoremstyle{remark}
\newtheorem{remark}[theorem]{Remark}
\newcommand{\Filp}{\mathsf{Fil}^{*}}
\newcommand{\Wpar}{\mathsf{W\text{-}Net}}
\newcommand{\Kpar}{\mathsf{K\text{-}Net}}
\newcommand{\Wex}{\mathsf{W\text{-}Net}_{\exists}}
\newcommand{\Kex}{\mathsf{K\text{-}Net}_{\exists}}
\newcommand{\NetAA}{\mathsf{AA\text{-}Net}}
\newcommand{\FF}{\mathcal F}
\newcommand{\GG}{\mathcal G}

\newcommand{\BB}{\mathcal B}
\newcommand{\CC}{\mathsf C}
\newcommand{\TT}{\mathcal T}
\newcommand{\NN}{\mathbb N}
\newcommand{\Ev}{\operatorname{Ev}}
\newcommand{\Tail}{\mathcal L}
\newcommand{\id}{\operatorname{id}}
\newcommand{\Fun}{\operatorname{Fun}}
\newcommand{\Iso}{\operatorname{Iso}}
\newcommand{\Sts}{S_{\mathrm{ts}}}
\newcommand{\Seq}{S_{\mathrm{eq}}}
\newcommand{\Ssat}{S_{\mathrm{sat}}}
\newcommand{\up}[1]{\mathord{\uparrow}#1}
\title{Filters as a localization of nets}
\author{Germ\'an I. Gomero Ferrer \and Renan M. Mezabarba}
\date{}
\hypersetup{pdftitle={Filters as a localization of nets},pdfauthor={German I. Gomero Ferrer and Renan M. Mezabarba}}
\begin{document}
\maketitle
\begin{abstract}
Based on different notions of subnet in the literature, we define
categories of nets and investigate their categorical relationships.
We show that the category of proper filters with filter-preserving
functions is, up to equivalence, a localization of four of these
categories.
\end{abstract}



\section*{Introduction}
This paper grew out of the authors' interest in convergence spaces,
both in research and in teaching, together with their enthusiasm for
category theory. Studying the so-called equivalence between nets and
filters in topological contexts, we asked how this relationship can be
expressed categorically, and how it depends on the notions of subnet
and morphism under consideration.

\emph{Nets} and \emph{filters} are classical tools for describing convergence beyond
sequences. Their relationship was studied by
Bruns and Schmidt~\cite{BS}, who gave net realizations of filters, and by
Aarnes and Anden\ae s~\cite{AA}, who compared notions of subnet through the
filters generated by \emph{tail sets}. Standard accounts appear in
G\"ahler~\cite{Gaehler}, Willard~\cite{Willard}, and
Schechter~\cite[Chapter 7]{Schechter}.

The convergence of a net in a topological space is determined by its
eventuality filter, which records the subsets containing all sufficiently
late terms of the net. Passing to this filter ``forgets'' the chosen
directed indexing of the net. This correspondence lets convergence
statements be translated between the two languages, but it does not by
itself say in what sense, and to what extent, the two languages are
equivalent. We want to know which information survives the passage
from nets to filters, and in what categorical sense this passage is
universal.
We extend the classical subnet relations to functions between possibly
different sets and investigate the resulting categories of nets.

We show that this relationship can be understood categorically through
localization. More precisely, the category of proper filters and
filter-preserving functions is, up to equivalence, a localization of
each of four categories of nets. This is stated in
Theorem~\ref{thm:fourlocal} and Corollary~\ref{cor:filterslocal}. In each case, the
localization can be obtained by formally inverting morphisms whose
value function is an identity and whose source and target have the
same eventuality filter. These terms are defined precisely in
Section~\ref{sec:categories}.

The paper is organized as follows.
Section~\ref{sec:categories} recalls the notions of nets, filters, and
subnets, and defines the categories to be compared.
Section~\ref{sec:canonical} treats realizations of filters by nets.
Section~\ref{sec:graphs} represents morphisms of nets and compares
their representations and compositions.
Section~\ref{sec:local} proves the localization results and examines
right fractions. In Section~\ref{sec:final}, we return to nets and filters
on a fixed set and discuss further aspects of their relationship.

\begin{samepage}
\paragraph{Acknowledgment of AI assistance.}
The categories of nets studied here were formulated by the authors.
They used GPT-6 Astra (OpenAI) and Claude
(Anthropic) to help organize and discuss calculations and drafts
handwritten by the authors since 2024. The interactions with these
systems helped the authors recognize the localization described here.
As usual, the authors take responsibility for the mathematical content
and the final presentation.\par
\end{samepage}

\section{Nets, filters and morphisms}\label{sec:categories}
We recall the conventions used throughout the paper, mostly adapting
Schechter's terminology, whose work~\cite[Chapter 7]{Schechter}
remains a general reference for the relation between nets, filters
and the different notions of subnet.
A \emph{directed set} is a nonempty preordered set in which any two
elements have a common upper bound. Its order is thus reflexive and
transitive, but need not be antisymmetric. For a directed set \(\mathbb{D}\)
and \(d\in\mathbb{D}\), write
\[
 \mathbb{D}_{\ge d}=\{d'\in \mathbb{D}:d\le d'\}.
\]

A subset \(C\subseteq \mathbb{D}\) is \emph{cofinal} if every element of \(\mathbb{D}\)
is below some element of \(C\). A map between preordered sets is
\emph{monotone} if it preserves their orders. We call it
\emph{monotone cofinal} when it is monotone and has cofinal range.
A \emph{net in a set \(X\)} consists of a directed set
\((\mathbb{D},\le)\) and a function \(x:\mathbb{D}\to X\). Both the
directed domain and the set \(X\) are part of its data. The set
\(x[\mathbb{D}_{\ge d}]\) is its \emph{tail set} at
\(d\). We denote the family of all such sets by
\[
 \Tail(x)=\{x[\mathbb{D}_{\ge d}]:d\in \mathbb{D}\}.
\]

A \emph{filter} on \(X\) is a nonempty family of subsets of \(X\)
that is upward closed and closed under finite intersections. It is
\emph{proper} if it does not contain the empty set. All filters in
this paper are proper. A \emph{filter base} is a nonempty family
\(\BB\) of nonempty subsets of \(X\) such that, for any
\(A,B\in\BB\), some \(C\in\BB\) satisfies \(C\subseteq A\cap B\).
The filter generated by \(\BB\) is
\[
 \langle\BB\rangle=\{A\subseteq X:\text{some }B\in\BB
                                      \text{ satisfies }B\subseteq A\}.
\]

We call \(\BB\) a base of a filter \(\FF\) on \(X\) when
\(\langle\BB\rangle=\FF\). Given a function \(f:X\to Y\), the image
filter of \(\FF\) is
\[
 f(\FF)=\{B\subseteq Y:f^{-1}(B)\in\FF\}.
\]

It is proper whenever \(\FF\) is proper. To relate nets and filters,
observe that the tail sets of a net form a filter base, which generates
the \emph{eventuality filter}
\[
 \Ev(x)=\{A\subseteq X:\text{there is }d\in\mathbb{D}
                  \text{ with }x[\mathbb{D}_{\ge d}]\subseteq A\}.
\]

The following identities relate eventuality filters to image filters
and describe image filters under composition and identities
(here \(g:Y\to Z\) is another function):
\begin{equation}\label{eq:image}
 \begin{aligned}
 \Ev(f\circ x)&=f(\Ev(x)),\\
 (g\circ f)(\FF)&=g(f(\FF)),\\
 \id_X(\FF)&=\FF.
 \end{aligned}
\end{equation}

These image filters allow us to compare filters on different sets.
The category \(\Filp\) has pairs \((X,\FF)\), where \(\FF\) is a
proper filter on \(X\), as objects. Its morphisms
\((X,\FF)\to(Y,\GG)\) are functions \(f:X\to Y\) satisfying
\(\GG\subseteq f(\FF)\). Identities and composition are those of
functions. This is the proper-filter part of Blass's category
\(\mathcal F\)~\cite[Section 1]{Blass}, or equivalently the category
\(F_0\) of Koubek--Reiterman~\cite[p. 19]{KR} restricted to proper
filters.

Filters on a fixed set are naturally compared by inclusion. For nets
in that set, comparing eventuality filters amounts to asking whether
each tail set of one net contains a tail set of the other. This gives
the notion of subnet considered by Aarnes and Anden\ae s~\cite[Definition 2.1 and Lemma 2.2]{AA}.
The subsequence relation between sequences also motivates subnet
relations defined by maps between index sets. We recall three notions
of subnet.

\begin{definition}\label{def:subnets}
Let \(x:\mathbb{D}\to X\) and \(y:\mathbb{E}\to X\) be nets in
the same set.
\begin{enumerate}[label=(\roman*)]
\item Following Willard~\cite[Definition 11.2]{Willard}, we call
\(x\) a \emph{W-subnet} of \(y\) if \(x=y\circ r\) for some
monotone cofinal map \(r:\mathbb{D}\to\mathbb{E}\).
\item Following Kelley's convention~\cite{Kelley}, as stated in
Schechter~\cite[Section 7.15(b)]{Schechter}, we call \(x\) a
\emph{K-subnet} of \(y\) if \(x=y\circ r\) for some
\emph{eventually cofinal} map \(r:\mathbb{D}\to\mathbb{E}\).
This means that, for every \(e_0\in\mathbb{E}\), there is
\(d_0\in\mathbb{D}\) such that
\begin{equation}\label{eq:ec}
 r[\mathbb{D}_{\ge d_0}]\subseteq\mathbb{E}_{\ge e_0}.
\end{equation}

\item Following Aarnes--Anden\ae s~\cite[Definition 2.1 and
Lemma 2.2]{AA}, we call \(x\) an \emph{AA-subnet} of \(y\) if
\(\Ev(y)\subseteq\Ev(x)\), equivalently if each tail set of \(y\)
contains a tail set of \(x\).\footnote{This relation is called
\emph{quasi-subnet} in~\cite[Section 2]{OTW}. See also
Schechter~\cite[Section 7.14]{Schechter} for further antecedents.}
\end{enumerate}
\end{definition}

For monotone maps, cofinal range is equivalent to eventual cofinality:
after an index maps above \(e_0\), every larger index does too.
Without monotonicity, cofinal range alone is insufficient. For example,
the map \(r:\NN\to\NN\) defined by \(r(2n)=n\) and \(r(2n+1)=0\),
where \(\NN\) is the set of nonnegative integers, has cofinal range
but fails~\eqref{eq:ec} at \(e_0=1\). A W-subnet is a K-subnet, and a K-subnet
is an AA-subnet. Counterexamples separating AA-subnets from
K-subnets are given in Aarnes--Anden\ae s~\cite[Examples 2.4--2.5]{AA}.
See also Schechter~\cite[Section 7.16]{Schechter} for comparisons
of all three notions.

\begin{definition}\label{def:tailequivalence}
Nets in one set are \emph{tail-equivalent}, written \(x\sim y\),
if \(\Ev(x)=\Ev(y)\).
\end{definition}

To extend subnet relations to functions between possibly different
sets, consider nets \(x:\mathbb{D}\to X\) and \(y:\mathbb{E}\to Y\),
and maps \(r:\mathbb{D}\to\mathbb{E}\) and \(f:X\to Y\) satisfying
\(f\circ x=y\circ r\), thus we have the commutative square:
\begin{equation}\label{eq:square}
 \begin{tikzcd}[row sep=large,column sep=large]
 \mathbb{D}\arrow[r,"r"]\arrow[d,"x"']&\mathbb{E}\arrow[d,"y"]\\
 X\arrow[r,"f"']&Y
 \end{tikzcd}
\end{equation}

We call \(f\) the \emph{value function} of the square and \(r\) its
\emph{index map}. The map \(r\) is a \emph{W-witness} for \(f\) when it is
monotone cofinal, and a \emph{K-witness} when it is eventually cofinal.
Thus an \emph{index witness} is an index map satisfying the commuting
equation and the specified order condition. Taking \(f=\id_X\)
recovers the corresponding subnet relation.

There are two ways to use these squares as morphisms: take the pair
\((r,f)\) as a morphism, or take only the value function \(f\) and
require that an index map \(r\) exist which makes the square commute
and satisfies the corresponding order condition. Combining these two
options with the two order conditions gives four categories, and a
fifth compares eventuality filters directly.
\begin{definition}\label{def:categories}
The following categories have all nets as objects.
For \(x:\mathbb{D}\to X\) and \(y:\mathbb{E}\to Y\), their morphisms \(x\to y\) are:
\begin{enumerate}[label=(\roman*)]
\item in \(\Wpar\), commutative squares \((r,f)\) with \(r\) monotone cofinal and
\(f\circ x=y\circ r\).
\item in \(\Kpar\), commutative squares \((r,f)\) with \(r\) eventually cofinal and
\(f\circ x=y\circ r\).
\item in \(\Wex\), functions \(f:X\to Y\) admitting a W-witness.
\item in \(\Kex\), functions \(f:X\to Y\) admitting a K-witness.
\item in \(\NetAA\), functions \(f:X\to Y\) such that
\(\Ev(y)\subseteq f(\Ev(x))\).
\end{enumerate}
\end{definition}

The prefixes W, K and AA refer to the notions of subnet in
Definition~\ref{def:subnets}. We call morphisms in the two
W-categories W-morphisms, those in the two K-categories K-morphisms,
and those of \(\NetAA\) AA-morphisms.

\medskip
Clearly, in \(\NetAA\), \(f\circ x\) is an AA-subnet
of \(y\), which in terms of tail sets translates into requiring that
for every \(e_0\in\mathbb{E}\) there is \(d_0\in\mathbb{D}\) such that
\begin{equation}\label{eq:aaquant}
 f[x[\mathbb{D}_{\ge d_0}]]\subseteq y[\mathbb{E}_{\ge e_0}].
\end{equation}

In \(\Wpar\) and \(\Kpar\), identities are
\((\id_{\mathbb{D}},\id_X)\), and composition is
\((s,g)\circ(r,f)=(s\circ r,g\circ f)\). This composition is well
defined, since composites of monotone cofinal maps are monotone
cofinal and composites of eventually cofinal maps are eventually
cofinal. In the other three categories, identities and composition
are those of the value functions.
There are full forgetful functors
\(P_W:\Wpar\to\Wex\) and \(P_K:\Kpar\to\Kex\),
which forget the index witness in each square.
The implications between the subnet conditions also give faithful
inclusions \(\Wpar\hookrightarrow\Kpar\),
\(\Wex\hookrightarrow\Kex\) and \(\Kex\hookrightarrow\NetAA\).

Let \(\CC\) be any of the categories in
Definition~\ref{def:categories}(i)--(iv). Composing the inclusions
and, when needed, the forgetful functor gives a functor
\(Q_{\CC}:\CC\to\NetAA\).
Finally, by taking eventuality filters we obtain the functor
\(T:\NetAA\to\Filp\), defined by
\(T(x:\mathbb{D}\to X)=(X,\Ev(x))\) on objects and \(T(f)=f\)
on morphisms. The definition of an AA-morphism makes \(T(f)\) a
filter morphism, and~\eqref{eq:image} ensures compatibility with
identities and composition. Together, these functors form the
commutative diagram
\[
\begin{tikzcd}[column sep=2.2em,row sep=2.5em]
\Wpar \arrow[r,hook] \arrow[d,"P_W"'] &
\Kpar \arrow[d,"P_K"] \\
\Wex \arrow[r,hook] &
\Kex \arrow[r,hook,"Q_{\Kex}"] & \NetAA \arrow[r,"T"] & \Filp.
\end{tikzcd}
\]

The functor \(T\) is fully faithful, since the corresponding
hom-sets consist of exactly the same functions. Every proper
filter is the eventuality filter of a net, as recalled in
Proposition~\ref{prop:base}, so \(T\) is essentially surjective
and hence an equivalence.

\begin{example}\label{ex:witness}
The functors \(P_W\) and \(P_K\) need not be faithful.
For \(c:\NN\to\{*\}\), the constant one-point net, the index maps
\(r(n)=n\) and \(s(n)=2n\) give distinct W-morphisms \(c\to c\) with
value function \(\id_{\{*\}}\). Their images under \(P_W\), and also
under \(P_K\), agree.
\end{example}

\begin{example}\label{ex:swap}
The inclusion \(\Wex\hookrightarrow\Kex\) need not be full, even
though only value functions are retained. Use the permutation of
consecutive pairs from Schechter~\cite[Section 7.16(b)]{Schechter}:
let \(y_n=n\) in \(\NN\) and \(x_n=r(n)\), where
\(r(2k)=2k+1\) and \(r(2k+1)=2k\). The map \(r\) is eventually
cofinal but not monotone. For the value function \(\id_{\NN}\),
injectivity of \(y\) forces this witness. Hence the W-to-K inclusions
are not full, including for the categories \(\Wex\) and \(\Kex\).
\end{example}

\Needspace{9\baselineskip}
\begin{example}\label{ex:head}
An AA-morphism need not have an index witness.
Let \(X=\{0,1\}\). For a nonempty \(A\subseteq X\), write
\(\up A=\{B\subseteq X:A\subseteq B\}\) for the principal filter
generated by \(A\). Put \(x=(1,0,0,\ldots)\) and
\(y=(0,0,0,\ldots)\). Their filters both equal \(\up{\{0\}}\), but
no map \(r\) satisfies \(x=y\circ r\), because of \(x_0\).
Thus \(\id_X:x\to y\) is an AA-morphism without any index witness.
The role of initial terms in this distinction already appears in
\cite[Example 2.4]{AA}. The example shows that the inclusion
\(\Kex\to\NetAA\) is not full.
\end{example}

\begin{example}\label{ex:parity}
A subnet can have a strictly finer eventuality filter.
For \(y_n=n\bmod2\) in \(X=\{0,1\}\), take \(r(n)=2n\).
The W-subnet \(x=y\circ r\) is constantly zero, while
\(\Ev(y)=\{X\}\) and \(\Ev(x)=\up{\{0\}}\).
This W-morphism has identity value function and refines eventuality
strictly. Inverting every
subnet would therefore discard information retained by \(\NetAA\).
\end{example}

\begin{remark}
\v{S}lapal~\cite[Definition 1.1 and Example 1.2(1)]{Slapal1996}
uses categorical data to define subnets by factorization through maps
between index sets. His category \(\mathsf{Dir}\) has directed sets as objects and
eventually cofinal maps as morphisms. This choice gives the K-subnet
relation recalled in Definition~\ref{def:subnets}(ii).
In later work, he considers, for a fixed set \(X\), the comma category
\((F\downarrow X)\), where \(F:\mathsf{Dir}\to\mathbf{Set}\) is the
forgetful functor. Its objects are nets in \(X\), and its morphisms
are index maps exhibiting subnets~\cite[Definition 2.1 and Example 2.2(2)]{Slapal2016}.
\end{remark}

\begin{remark}\label{rem:filterarrows}
For a directed set \(\mathbb D\), write \(\TT_{\mathbb D}=\Ev(\id_{\mathbb D})\)
for the filter generated by its tail sets. As observed in Schechter~\cite[Section 7.15(b)]{Schechter}, the index
condition for Kelley morphisms can itself be expressed in terms of filters. For any map \(r:\mathbb{D}\to\mathbb{E}\),
\(\Ev(r)=r(\TT_{\mathbb{D}})\), and
\[
 r\text{ is eventually cofinal}\quad\Leftrightarrow\quad
 \TT_{\mathbb{E}}\subseteq r(\TT_{\mathbb{D}}).
\]

In fact, the preimage of each tail set of \(\mathbb{E}\) contains a tail set
of \(\mathbb{D}\) exactly when~\eqref{eq:ec} holds.
Consequently \(\Kpar\) is fully faithfully represented in the arrow category of \(\Filp\). The representation sends a net \(x\) to
the filter morphism
\[
 x:(\mathbb{D},\TT_{\mathbb{D}})\to(X,\Ev(x))
\]
and a Kelley morphism \((r,f)\) to its commuting square, as
in~\eqref{eq:square}. Conversely,
for nets \(x:\mathbb D\to X\) and \(y:\mathbb E\to Y\), a square
between their associated objects in the arrow category has a
component \(r:(\mathbb D,\TT_{\mathbb D})\to
(\mathbb E,\TT_{\mathbb E})\). Since this component is a filter
morphism, the equivalence above makes \(r\) eventually cofinal.
The square therefore comes from a unique Kelley morphism.
If \(J\) denotes this representation and \(\operatorname{cod}\)
the codomain functor, then \(\operatorname{cod}\circ J=T\circ Q_{\Kpar}\).
\end{remark}

\section{Realizations of filters}\label{sec:canonical}
To compare the categories of nets with \(\Filp\), we need to realize
proper filters as eventuality filters of nets. We use the construction of Bruns--Schmidt~\cite[Section 3,
pp. 170--171]{BS} in its formulation for filter bases. See also
Schechter~\cite[Section 7.11]{Schechter},
Di Biase--Krantz~\cite[Lemma 6.3]{DBK}, and
O'Brien--Troitsky--van der Walt~\cite[Proposition 5.1]{OTW}.

For a filter base \(\BB\) on a set \(X\), let
\(\mathbb{D}_{\BB}=\{(A,a):A\in\BB,\ a\in A\}\), ordered by
\((A,a)\le(B,b)\) when \(B\subseteq A\). Since a proper filter
\(\FF\) is a base for itself, taking \(\BB=\FF\) gives the
construction of Bruns and Schmidt for \(\FF\).

\begin{proposition}\label{prop:base}
The pair \((\mathbb{D}_{\BB},\le)\) is a directed preorder. The net \(\nu_{\BB}:\mathbb{D}_{\BB}\to X\)
defined by \(\nu_{\BB}(A,a)=a\) satisfies
\begin{equation}\label{eq:exacttail}
 \nu_{\BB}[(\mathbb{D}_{\BB})_{\ge(A,a)}]=A
\end{equation}
for every \((A,a)\in\mathbb{D}_{\BB}\). Consequently,
\(\Tail(\nu_{\BB})=\BB\) and
\(\Ev(\nu_{\BB})=\langle\BB\rangle\).
\end{proposition}
\begin{proof}
For \((A,a),(B,b)\), choose \(C\in\BB\) and \(c\in C\)
in such a way that \(C\subseteq A\cap B\). Then \((C,c)\)
is a common upper bound. For every pair \((B,b)\ge(A,a)\),
we have \(B\subseteq A\), and hence
\(\nu_{\BB}(B,b)=b\in B\subseteq A\). Conversely,
\((A,a')\ge(A,a)\) for each \(a'\in A\), proving~\eqref{eq:exacttail}.
The assertions about tail sets and the generated filter follow.
\end{proof}

Every net with eventuality filter \(\FF\) is a W-subnet of
\(\nu_{\FF}\).
This is the equal-filter case of the comparison in
Costantini--Priola~\cite[Proposition 3.1]{CP}, whose subnet relation
is the Kelley relation. Their index map is in fact monotone.
We record this map together with the comparison obtained from a
filter base.

\begin{proposition}\label{prop:basecomparison}
Let \(\BB\) be a base of a proper filter \(\FF\) on \(X\).
The inclusion \(j:\mathbb{D}_{\BB}\to\mathbb{D}_{\FF}\) defines
a W-morphism \((j,\id_X):\nu_{\BB}\to\nu_{\FF}\) between nets
with the same eventuality filter. For every net
\(x:\mathbb{D}\to X\), the map
\(\mu_x:\mathbb{D}\to\mathbb{D}_{\Tail(x)}\) defined by
\[
 \mu_x(d)=(x[\mathbb{D}_{\ge d}],x_d)
\]
is monotone cofinal and satisfies
\(\nu_{\Tail(x)}\circ\mu_x=x\).
\end{proposition}
\begin{proof}
For \((A,a)\in \mathbb{D}_{\FF}\), choose \(B\in\BB\) with \(B\subseteq A\)
and \(b\in B\). Then \(j(B,b)\ge(A,a)\), proving cofinality.
Monotonicity follows from the definition of the orders. The commuting
equation and equality of filters follow from
Proposition~\ref{prop:base}.
For \(\mu_x\), inclusion between tail sets gives monotonicity. Given
\((x[\mathbb{D}_{\ge d_0}],a)\), the index \(d_0\) itself maps above it.
The commuting equation follows from the definition of \(\nu_{\Tail(x)}\),
and equality of filters follows from Proposition~\ref{prop:base} and
\(\langle\Tail(x)\rangle=\Ev(x)\).
\end{proof}

Notice that in the previous proof, \(\mu_x\) may be followed by
the inclusion \(j:\mathbb D_{\Tail(x)}\to\mathbb D_{\Ev(x)}\).
The result is a W-morphism \((j\circ\mu_x,\id_X):x\to\nu_{\Ev(x)}\)
with identity value function between two nets with eventuality
filter \(\Ev(x)\).

The nets \(\nu_{\FF}\) give an explicit quasi-inverse
\(N:\Filp\to\NetAA\), defined by \(N(X,\FF)=\nu_{\FF}\) on
objects and \(N(f)=f\) on morphisms.
Proposition~\ref{prop:base} gives \(T\circ N=\id_{\Filp}\), since
\(\Ev(\nu_{\FF})=\FF\) for every proper filter \(\FF\), and both
functors retain the function representing each morphism.
On the other hand, for each net \(x:\mathbb D\to X\), the identity
\(\id_X\) gives mutually inverse AA-morphisms between
\(\nu_{\Ev(x)}\) and \(x\). Denote the morphism from
\(\nu_{\Ev(x)}\) to \(x\) by \(\varepsilon_x\). For every
AA-morphism \(f:x\to y\), the diagram
\[
\begin{tikzcd}[row sep=large,column sep=large]
N(T(x)) \arrow[r,"N(T(f))"] \arrow[d,"\varepsilon_x"'] &
N(T(y)) \arrow[d,"\varepsilon_y"] \\
x \arrow[r,"f"'] & y
\end{tikzcd}
\]
commutes because both composites are represented by the same
function \(f\), that is,
\[
 \varepsilon_y\circ N(T(f))=f\circ\varepsilon_x.
\]

Thus these components form a natural isomorphism
\(\varepsilon:N\circ T\Rightarrow\id_{\NetAA}\).

\Needspace{18\baselineskip}
\section{``Localizing'' the missing ``arrow''}\label{sec:graphs}
The functor \(T:\NetAA\to\Filp\) is an equivalence. However, while in
a morphism of \(\Filp\)
\[
\begin{tikzcd}[column sep=large]
(X,\FF) \arrow[r,"f"] & (Y,\GG)
\end{tikzcd}
\]
everything looks fine, in an AA-morphism \(f:x\to y\)
\[
\begin{tikzcd}[column sep=large,row sep=large]
\mathbb D \arrow[r,dotted,"?"] \arrow[d,"x"'] & \mathbb E \arrow[d,"y"] \\
X \arrow[r,"f"'] & Y
\end{tikzcd}
\]
the picture looks incomplete: an index map, or something playing its
role, is missing.

We show that an index map becomes available
after replacing the source by a tail-equivalent W-subnet. The
construction extends the one in Schechter~\cite[Lemma 7.18 and
Corollary 7.19]{Schechter} and
O'Brien--Troitsky--van der Walt~\cite[Proposition 2.2]{OTW}
to a function between the sets in which the nets take values.\footnote{For
a related realization of a finer filter by a subnet, see
G\"ahler~\cite[Theorem 2.8.4, pp. 82--83]{Gaehler}.}

An eventually cofinal map \(r:\mathbb A\to\mathbb D\) between
directed sets makes \(x\circ r\) a K-subnet of \(x\) for every
net \(x:\mathbb D\to X\). Thus \(\Ev(x)\subseteq\Ev(x\circ r)\).
To obtain \(\Ev(x\circ r)=\Ev(x)\) for every such net, we also
need the reverse inclusion between the index filters. This motivates
the following condition.

\begin{definition}\label{def:tailsurj}
We say that a monotone map \(r:\mathbb A\to\mathbb D\) between directed sets is
\emph{tail-surjective} if
\begin{equation}\label{eq:tailsurjfilters}
 r(\TT_{\mathbb A})\subseteq\TT_{\mathbb D}.
\end{equation}

Equivalently, the image of each tail set of \(\mathbb A\) contains a tail set
of \(\mathbb D\): for each \(a_0\in\mathbb A\), there is
\(d_0\in\mathbb D\) such that
\(\mathbb D_{\ge d_0}\subseteq r[\mathbb A_{\ge a_0}]\).
\end{definition}

The image of a tail-surjective map is cofinal in its codomain.
Indeed, fix \(a_0\), choose
\(d_0\) as above, and let \(d\in\mathbb D\). An upper bound of
\(d,d_0\) belongs to \(r[\mathbb A_{\ge a_0}]\). Since the map is
monotone, it is therefore eventually cofinal as well. Consequently,
\eqref{eq:tailsurjfilters} is an equality for a tail-surjective map.

\begin{proposition}\label{prop:binary}
For a monotone map \(r:\mathbb A\to\mathbb D\), the following
conditions are equivalent:
\begin{enumerate}[label=(\roman*)]
\item \(r\) is tail-surjective;
\item \(\Ev(x\circ r)=\Ev(x)\) for every net \(x:\mathbb D\to X\),
on every set \(X\);
\item \(\Ev(x\circ r)=\Ev(x)\) for every net
\(x:\mathbb D\to\{0,1\}\).
\end{enumerate}
\end{proposition}
\begin{proof}
Suppose (i) holds. Eventual cofinality gives
\(\Ev(x)\subseteq\Ev(x\circ r)\). If \(x\circ r\) takes its values
in \(B\) after an index \(a_0\), choose \(d_0\) such that every
\(d\ge d_0\) has the form \(r(a)\) for some \(a\ge a_0\).
Then \(x_d\in B\) for every \(d\ge d_0\). This proves the reverse
inclusion and hence (ii), which implies (iii).

For the remaining implication, suppose tail-surjectivity fails at
\(a_0\). Then
\(H=\mathbb D\setminus r[\mathbb A_{\ge a_0}]\) meets every tail set
of \(\mathbb D\). Its characteristic net
\(x=\mathbf1_H:\mathbb D\to\{0,1\}\) is not eventually zero,
whereas \(x\circ r\) is zero after \(a_0\). Their eventuality
filters differ, contradicting (iii).
\end{proof}

To supply the missing index map, we will replace a net
\(x:\mathbb D\to X\) by a net \(x\circ\sigma\) with the same
eventuality filter, where \(\sigma:\mathbb A\to\mathbb D\) is a
monotone map between directed sets. By Proposition~\ref{prop:binary},
this holds for every net indexed by \(\mathbb D\) exactly when
\(\sigma\) is tail-surjective.

For nets \(x:\mathbb D\to X\) and \(y:\mathbb E\to Y\), we say that a morphism
\(x\to y\) in any of the categories of Definition~\ref{def:categories}
is \emph{vertical} if \(X=Y\) and its value function is \(\id_X\).
A vertical morphism is \emph{exact} if \(x\sim y\). We then write
\begin{equation}\label{eq:sts}
 \Sts=\{(r,\id_X):r\text{ is tail-surjective}\},
\end{equation}
where \((r,\id_X)\) ranges over vertical morphisms of \(\Wpar\).
Every member of \(\Sts\) is exact by
Proposition~\ref{prop:binary}. The letter \(S\) anticipates the
notation \(\CC[S^{-1}]\) for localization at a class of morphisms:
\(\Sts\) will be one such class, with the subscript recording
tail-surjectivity. If necessary, we may regard \(\Sts\) as a class
of morphisms in \(\Kpar\) through the inclusion
\(\Wpar\hookrightarrow\Kpar\).

The class \(\Sts\) contains identities and is closed under
composition. In fact, if \(r:\mathbb A\to\mathbb D\) and
\(s:\mathbb D\to\mathbb E\) are tail-surjective, then for each
\(a_0\) there are \(d_0,e_0\) with
\[
 \mathbb E_{\ge e_0}\subseteq s[\mathbb D_{\ge d_0}]
                  \subseteq(s\circ r)[\mathbb A_{\ge a_0}].
\]

This property will allow successive replacements of index domains.
We now use these morphisms to compare an AA-morphism \(f:x\to y\)
with morphisms in the categories \(\Wpar\) and \(\Kpar\). The comparison will have
a common source \(u\), with a morphism \(s:u\to x\) in \(\Sts\) and
a morphism \(w:u\to y\) whose value function is \(f\).
Thus \(u\) replaces \(x\) without changing its eventuality filter,
and \(w\) has an index map compatible with \(f\).\footnote{Such diagrams
also occur in localization theory. See Mili\v{c}i\'c~\cite[Chapter 1, Section 1.3]{Milicic}.}

\begin{definition}\label{def:cones}
For nets \(x:\mathbb D\to X\) and \(y:\mathbb E\to Y\), fix an
AA-morphism \(f:x\to y\) and let \(\CC\) be \(\Wpar\) or \(\Kpar\).
A \emph{cone representing \(f\) in \(\CC\)} is a diagram in \(\CC\)
\begin{equation}\label{eq:representingcone}
\begin{tikzcd}[column sep=2em,row sep=1.5em]
& u \arrow[dl,"s"'] \arrow[dr,"w"] & \\
x & & y
\end{tikzcd}
\end{equation}
where \(u\) is a net in \(X\), \(s\in\Sts\) and the value function
of \(w\) is \(f\). We denote by \(\CC^{\triangle f}\) the
category whose objects are these cones and whose morphisms from
\((u,s,w)\) to \((u',s',w')\) are vertical \(\CC\)-morphisms
\(h:u\to u'\) satisfying
\[
 s'\circ h=s \qquad\text{and}\qquad w'\circ h=w.
\]

Identities and composition are inherited from \(\CC\).
\end{definition}

For a cone \((u,s,w)\) as in~\eqref{eq:representingcone}, write
\(u:\mathbb A\to X\), \(s=(\sigma,\id_X)\) and \(w=(k,f)\). Here
\(\sigma:\mathbb A\to\mathbb D\) is tail-surjective, and
\(k:\mathbb A\to\mathbb E\) is monotone cofinal when \(\CC=\Wpar\) and
eventually cofinal when \(\CC=\Kpar\). The cone thus replaces the
incomplete square of \(f\) from the beginning of this section by two
commuting squares:
\begin{equation}\label{eq:conecomparison}
\begin{tikzcd}[column sep=large,row sep=large,baseline={([yshift=-.5ex]current bounding box.center)}]
\mathbb D \arrow[r,dotted,"?"] \arrow[d,"x"'] & \mathbb E \arrow[d,"y"] \\
X \arrow[r,"f"'] & Y
\end{tikzcd}
\qquad\rightsquigarrow\qquad
\begin{tikzcd}[column sep=large,row sep=large,baseline={([yshift=-.5ex]current bounding box.center)}]
\mathbb D \arrow[d,"x"'] &
\mathbb A \arrow[l,"\sigma"'] \arrow[r,"k"] \arrow[d,"u"] &
\mathbb E \arrow[d,"y"] \\
X & X \arrow[l,"\id_X"] \arrow[r,"f"'] & Y
\end{tikzcd}
\end{equation}

The dotted side in~\eqref{eq:conecomparison} marks the missing map. The square with sides \(\sigma\), \(x\), \(u\) and
\(\id_X\) is the morphism \(s\), and the square with sides \(k\),
\(u\), \(y\) and \(f\) is the morphism \(w\). Thus \(u=x\circ\sigma\),
which has the same eventuality filter as \(x\), and
\(f\circ u=y\circ k\). Although \(f\) need not have an index witness,
the pair \(\sigma,k\) plays the role of the missing index map after
\(x\) is replaced by the tail-equivalent net \(u\). In
Section~\ref{sec:local}, once \(s\) is inverted, the cone yields the
morphism \(w\circ s^{-1}\). This is the sense in which the missing
``arrow'' is localized.

\begin{theorem}\label{thm:graph}\label{prop:terminalroof}
For an AA-morphism \(f:x\to y\) and either \(\CC=\Wpar\)
or \(\CC=\Kpar\), the category \(\CC^{\triangle f}\) has a
terminal object. More precisely, there is a representing cone
\begin{equation}\label{eq:cone}
\begin{tikzcd}[column sep=2em,row sep=1.5em]
& z_f \arrow[dl,"s_f"'] \arrow[dr,"w_f"] & \\
x & & y
\end{tikzcd}
\end{equation}
in which \(s_f\) and \(w_f\) are W-morphisms, and which is
terminal for both choices of \(\CC\). In particular,
\(\Ev(z_f)=\Ev(x)\).
\end{theorem}
\Needspace{9\baselineskip}
\begin{proof}
Write \(x:\mathbb D\to X\) and \(y:\mathbb E\to Y\). Give
\[
 C_f=\{(d,e)\in\mathbb D\times\mathbb E:f(x_d)=y_e\}
\]
the induced product order. We first show that \(C_f\) is cofinal
in \(\mathbb D\times\mathbb E\). Given \((d_0,e_0)\) in this product,
the AA-condition gives an index \(d_*\) such that, for every
\(d\ge d_*\), there is an index \(e\ge e_0\) with \(f(x_d)=y_e\).
Choose \(\tilde d\ge d_0,d_*\) and then \(e\ge e_0\) with
\(f(x_{\tilde d})=y_e\). Then \((\tilde d,e)\in C_f\) lies above
\((d_0,e_0)\).

As a cofinal subset of a directed set, \(C_f\) is directed. Its
coordinate projections \(p:C_f\to\mathbb D\) and
\(q:C_f\to\mathbb E\) are monotone cofinal. To finish the setup, we
define the net \(z_f:C_f\to X\) by \(z_f(d,e)=x_d\) for every
\((d,e)\in C_f\). These maps form the commutative diagram
\begin{equation}\label{eq:representationdiagram}
\begin{tikzcd}[column sep=large,row sep=large]
\mathbb D \arrow[r,"x"] & X \arrow[dd,"f"] \\
C_f \arrow[u,"p"] \arrow[d,"q"']
    \arrow[ur,"z_f" description]
    \arrow[dr,"f\circ z_f" description] & \\
\mathbb E \arrow[r,"y"'] & Y.
\end{tikzcd}
\end{equation}

Consequently \(s_f=(p,\id_X)\) and \(w_f=(q,f)\) are W-morphisms.
The projection \(p\) is tail-surjective. Indeed, for
\(c_0=(d_0,e_0)\in C_f\), choose \(d_*\) as above and then
\(d_1\ge d_0,d_*\). Every \(d\ge d_1\) admits some \(e\ge e_0\)
with \(f(x_d)=y_e\), and then \((d,e)\in(C_f)_{\ge c_0}\), so
\begin{equation}\label{eq:stronggraph}
 \mathbb D_{\ge d_1}\subseteq p[(C_f)_{\ge c_0}].
\end{equation}

This proves \(s_f\in\Sts\). Proposition~\ref{prop:binary} now
gives \(\Ev(z_f)=\Ev(x)\).
To prove terminality, consider any cone \((u,s,w)\) representing
\(f\), and write \(u:\mathbb A\to X\),
\(s=(\sigma,\id_X)\), and \(w=(k,f)\). The square equations
give \(f(x_{\sigma(a)})=f(u_a)=y_{k(a)}\) for every
\(a\in\mathbb A\). Hence the map
\(\langle\sigma,k\rangle:\mathbb A\to C_f\) defined by
\[
 \langle\sigma,k\rangle(a)=(\sigma(a),k(a))
\]
is well defined. Both coordinates are eventually cofinal. For a
given \((d_0,e_0)\in C_f\), they therefore lie above
\(d_0,e_0\), respectively, on a common tail set of \(\mathbb A\).
Thus \(\langle\sigma,k\rangle\) is eventually cofinal. In \(\Wpar\), both
coordinates are also monotone, so \(\langle\sigma,k\rangle\) is monotone cofinal.
It follows that \(h=(\langle\sigma,k\rangle,\id_X):u\to z_f\) is a
\(\CC\)-morphism, and its coordinate projections give
\(s_f\circ h=s\) and \(w_f\circ h=w\). These equations
uniquely determine both components of \(h\), as desired.
\end{proof}

By Theorem~\ref{thm:graph}, every cone \((u,s,w)\) representing \(f\)
has a unique morphism to \((z_f,s_f,w_f)\) in \(\CC^{\triangle f}\).
Its underlying \(\CC\)-morphism \(h:u\to z_f\) satisfies
\(s=s_f\circ h\) and \(w=w_f\circ h\), so all these cones determine
the same morphism once the morphisms in \(\Sts\) are sent to
isomorphisms. Indeed, let \(\Phi:\CC\to\mathsf A\) be a functor that
sends \(\Sts\) to isomorphisms. For every cone \((u,s,w)\)
representing \(f\),
\begin{equation}\label{eq:roofindependence}
 \Phi(w)\circ\Phi(s)^{-1}
       =\Phi(w_f)\circ\Phi(s_f)^{-1}.
\end{equation}

To see this, note that \(\Phi(s_f)\circ\Phi(h)=\Phi(s)\), and both
\(\Phi(s_f)\) and \(\Phi(s)\) are invertible. Thus
\(\Phi(h)\) is invertible too. Substitution of
\(s=s_f\circ h\) and \(w=w_f\circ h\) proves
\eqref{eq:roofindependence}.

\begin{samepage}
\begin{corollary}\label{cor:replacement}
For fixed nets \(x:\mathbb D\to X\) and \(y:\mathbb E\to Y\), and a
function \(f:X\to Y\), the following are equivalent:
\begin{enumerate}[label=(\roman*)]
\item \(f:x\to y\) is an AA-morphism;
\item there are a net \(x'\) with \(x'\sim x\) and a W-morphism
\((r,f):x'\to y\);
\item there are nets \(x'\) and \(y'\) with \(x'\sim x\) and
\(y'\sim y\), and a W-morphism \(x'\to y'\) with value function \(f\).
\end{enumerate}
In (ii), \(x'\) can be chosen with a vertical W-morphism
\(x'\to x\) belonging to \(\Sts\).
\end{corollary}
\end{samepage}
\begin{proof}
Theorem~\ref{thm:graph} proves (i) implies (ii), and (ii) implies
(iii) by taking \(y'=y\). A W-morphism as in (iii) gives
\(\Ev(y')\subseteq f(\Ev(x'))\). Replacing these filters by
\(\Ev(y)\) and \(\Ev(x)\), respectively, proves (i).
\end{proof}

We next relate representations to composition. Let \(f:x\to y\) and
\(g:y\to v\) be AA-morphisms. Theorem~\ref{thm:graph} gives terminal
representing cones \((z_f,s_f,w_f)\) and \((z_g,s_g,w_g)\) in
\(\CC^{\triangle f}\) and \(\CC^{\triangle g}\), respectively.
From these we shall construct a cone representing \(g\circ f\).
To build the new cone, we complete a square with sides \(w_f:z_f\to y\)
and \(s_g:z_g\to y\), obtaining the required common source.

\begin{proposition}\label{prop:ore}
Let \(\CC\) be \(\Wpar\) or \(\Kpar\), and suppose
\(b=(r,f):x\to y\) and \(s=(\sigma,\id_Y):z\to y\) are
\(\CC\)-morphisms with \(s\in\Sts\). There exist a W-morphism
\(h:u\to z\) with value function \(f\) and a morphism
\(t:u\to x\) in \(\Sts\) making the following square commute
in \(\CC\):
\[
\begin{tikzcd}[column sep=large,row sep=large]
u \arrow[r,"h"] \arrow[d,"t"'] & z \arrow[d,"s"] \\
x \arrow[r,"b"'] & y.
\end{tikzcd}
\]
\end{proposition}
\begin{proof}
Let \(\mathbb D,\mathbb E,\mathbb A\) be the index sets of
\(x,y,z\), respectively. Give
\[
 G=\{(d,a)\in\mathbb D\times\mathbb A:r(d)=\sigma(a)\}
\]
the induced product order. Fix \((d_0,a_0)\) in the product.
By tail-surjectivity of \(\sigma\), there is an index \(e_1\) with
\(\mathbb E_{\ge e_1}\subseteq\sigma[\mathbb A_{\ge a_0}]\).
Eventual cofinality of \(r\) then gives an index \(d_1\ge d_0\)
such that \(r(d)\ge e_1\) for every \(d\ge d_1\).
For each such \(d\), there is \(a\ge a_0\) with
\(r(d)=\sigma(a)\). Thus \(G\)
is cofinal in the product and hence directed. Its projections
\(p:G\to\mathbb D\) and \(q:G\to\mathbb A\) are monotone
cofinal. Applied to \((d_0,a_0)\in G\), the same argument gives
\[
 \mathbb D_{\ge d_1}\subseteq p[G_{\ge(d_0,a_0)}],
\]
so \(p\) is tail-surjective.

Define \(u(d,a)=x_d\). On \(G\),
\[
 f(x_d)=y_{r(d)}=y_{\sigma(a)}=z_a.
\]

Consequently \(t=(p,\id_X)\) and \(h=(q,f)\) are W-morphisms,
with \(t\in\Sts\). The index maps of their composites satisfy
\(r\circ p=\sigma\circ q\), and their value functions satisfy
\(f\circ\id_X=\id_Y\circ f\). This proves
\(b\circ t=s\circ h\) in either category.
\end{proof}

For the cones of Theorem~\ref{thm:graph}, apply
Proposition~\ref{prop:ore} to \(b=w_f\) and \(s=s_g\).
The resulting cone
\[
\begin{tikzcd}[column sep=2.5em,row sep=1.8em]
& u \arrow[dl,"s_f\circ t"'] \arrow[dr,"w_g\circ h"] & \\
x & & v
\end{tikzcd}
\]
represents \(g\circ f\), since \(s_f\circ t\in\Sts\).
By terminality in Theorem~\ref{thm:graph}, this cone has a unique
morphism to the chosen cone representing \(g\circ f\). The
corresponding instance of~\eqref{eq:roofindependence} is the identity
used for composition in the proof of Theorem~\ref{thm:wlocal}.

\begin{remark}
For either choice of \(\CC\) in Proposition~\ref{prop:ore}, the square
constructed in its proof is a pullback in \(\CC\). To see this, let
\(\xi:\mathbb B\to Z\)
be a net with morphisms \(t'=(\alpha,\varphi):\xi\to x\) and
\(h'=(\beta,\psi):\xi\to z\) in the chosen category, and suppose
\(b\circ t'=s\circ h'\). Equality of the two components gives
\(r\circ\alpha=\sigma\circ\beta\) and \(\psi=f\circ\varphi\).
The map \(\ell:\mathbb B\to G\), defined by
\(\ell(i)=(\alpha(i),\beta(i))\), is therefore well defined.
Since \(\alpha\) and \(\beta\) are eventually cofinal, for every
\((d,a)\in G\) both coordinates of \(\ell\) are eventually above
\(d\) and \(a\), respectively. Hence \(\ell\) is eventually cofinal.
When the two maps are monotone, so is \(\ell\). Moreover,
\(u\circ\ell=x\circ\alpha=\varphi\circ\xi\), so
\((\ell,\varphi):\xi\to u\) is a morphism in the chosen category.
Its composites with \(t\) and \(h\) are \(t'\) and \(h'\), respectively.
The equations with the projections determine both components of this
morphism, so it is unique. Thus the pullback constructed here has
\(t\in\Sts\). This square establishes for \(\Sts\) the condition (F2)
introduced in Section~\ref{sec:local}.
\end{remark}

\section{Localization, in the categorical sense}\label{sec:local}
We use the universal formulation of localization in
Gabriel--Zisman~\cite[Chapter I, Section 1.1 and Lemma 1.2]{GZ}.
If \(S\) is a class of morphisms in \(\CC\), a functor
\(q:\CC\to\mathsf L\) is a \emph{localization at \(S\)} if
it sends \(S\) to isomorphisms and, for every category
\(\mathsf A\), precomposition gives an equivalence
\begin{equation}\label{eq:universal}
 q^*:\Fun(\mathsf L,\mathsf A)\to\Fun_S(\CC,\mathsf A).
\end{equation}

Here \(\Fun_S(\CC,\mathsf A)\) is the full subcategory of
functors sending \(S\) to isomorphisms. Morphisms in a class being
inverted are called \emph{denominators}. For the functors
\(Q_{\CC}\) considered below, precomposition will in
fact be an isomorphism of functor categories.

We begin with the category \(\Wpar\). Write
\(Q=Q_{\Wpar}:\Wpar\to\NetAA\) for the functor that
sends a commutative square \((r,f)\) to its value
function \(f\) (cf. Section~\ref{sec:categories}).
Previously, we introduced \(\Sts\) as the class of vertical morphisms
whose index maps satisfy the tail-surjectivity condition used in
Theorem~\ref{thm:graph} and Proposition~\ref{prop:ore}.
We now compare it with a larger class. For nets \(x:\mathbb D\to X\) and
\(y:\mathbb E\to X\), a vertical morphism \((r,\id_X):x\to y\) is a
commutative square
\begin{equation}\label{eq:verticalsquare}
 \begin{tikzcd}[row sep=large,column sep=large]
 \mathbb{D}\arrow[r,"r"]\arrow[d,"x"']&\mathbb{E}\arrow[d,"y"]\\
 X\arrow[r,"\id_X"']&X
 \end{tikzcd}
\end{equation}
and preserving eventuality requires only
\(\Ev(x)=\Ev(y)\). So we write \(\Seq\) for this class of
exact vertical morphisms in \(\Wpar\), with the subscript
recording equality of the two eventuality filters:
\begin{equation}\label{eq:seq}
 \Seq=\{(r,\id_X):x\sim y\}.
\end{equation}

Proposition~\ref{prop:binary} gives \(\Sts\subseteq\Seq\).
Since \(x=y\circ r\), membership in \(\Seq\) requires
\(\Ev(y\circ r)=\Ev(y)\) only for the net \(y\). By
Proposition~\ref{prop:binary}, tail-surjectivity of \(r\) requires
\(\Ev(y'\circ r)=\Ev(y')\) for every net \(y'\) indexed by
\(\mathbb E\).

Any class of denominators for which \(Q\) is a localization must
consist of morphisms sent to isomorphisms by \(Q\). The largest possible
such class is therefore the preimage of the isomorphisms of
\(\NetAA\). We first describe these isomorphisms.

\begin{lemma}\label{lem:iso}
An AA-morphism \(f:x\to y\) is an isomorphism precisely when
\(f\) is bijective and \(f(\Ev(x))=\Ev(y)\).
\end{lemma}
\begin{proof}
An inverse morphism must be an inverse function, so \(f\) is bijective.
The morphism conditions for \(f\) and \(f^{-1}\) then give equality
of the filters. Conversely, the stated equality makes \(f^{-1}\)
an AA-morphism, and the two composites are identity functions.
\end{proof}

We therefore put
\begin{equation}\label{eq:ssat}
 \Ssat=Q^{-1}(\Iso(\NetAA)).
\end{equation}

Combining the definitions with Lemma~\ref{lem:iso}, which characterizes
the isomorphisms of \(\NetAA\), gives
\begin{equation}\label{eq:classes}
 \Sts\subseteq\Seq\subseteq\Ssat.
\end{equation}

Each class contains identities and is closed under composition.
For morphisms \(x\to y\) and \(y\to z\) both in \(\Seq\),
we have \(\Ev(x)=\Ev(y)=\Ev(z)\), so their composite also
belongs to \(\Seq\). For \(\Ssat\), closure follows because
isomorphisms are closed under composition.

We shall prove that \(\Ssat\) is the \emph{saturation} of
\(\Sts\) and also of \(\Seq\) in the following sense: it consists of
the morphisms that are sent to isomorphisms by every functor inverting
\(\Sts\), respectively \(\Seq\). These morphisms are forced to become invertible,
and adding them does not change the localization. This will follow
from the universal property proved below.

\begin{example}\label{ex:classes}
For the one-point constant net indexed by \(\NN\), the index map
\(n\mapsto2n\) defines an exact vertical W-morphism. It is not
tail-surjective, because its image misses every odd index. Hence
\(\Sts\subsetneq\Seq\). A bijection between two distinct
singleton sets, equipped with singleton-indexed nets, belongs to
\(\Ssat\) and is not vertical. Thus \(\Seq\subsetneq\Ssat\).

Proposition~\ref{prop:basecomparison} gives another exact morphism
that need not belong to \(\Sts\). Let \(\FF\) be the cofinite
filter on \(\NN\), and let
\(\BB=\{\{k\in\NN:k\ge n\}:n\in\NN\}\).
The base inclusion \(\mathbb D_{\BB}\to\mathbb D_{\FF}\)
is exact. Every tail set of \(\mathbb D_{\FF}\), however, contains
a pair \((A,k)\) in which \(A\) is a sufficiently late numerical
tail set with one of its noninitial points removed. Such an \(A\)
does not belong to \(\BB\), so the inclusion is not
tail-surjective.
\end{example}

\begin{theorem}\label{thm:wlocal}
The functor \(Q:\Wpar\to\NetAA\) is a localization at each of the
classes \(\Sts\), \(\Seq\) and \(\Ssat\). More precisely, for every
category \(\mathsf A\), the functor in \eqref{eq:universal} is an
isomorphism.
\end{theorem}
\begin{proof}
By \eqref{eq:classes}, \(Q\) sends each of the three classes to
isomorphisms. Suppose first that \(\Phi:\Wpar\to\mathsf A\) inverts
\(\Sts\). Define \(\overline\Phi\) on objects by
\(\overline\Phi(x)=\Phi(x)\) and on morphisms by
\begin{equation}\label{eq:extension}
 \overline\Phi(f)=\Phi(w_f)\circ\Phi(s_f)^{-1},
\end{equation}
where \((z_f,s_f,w_f)\) is the cone of Theorem~\ref{thm:graph}.
By~\eqref{eq:roofindependence}, every cone \((u,s,w)\) representing
\(f\) satisfies \(\Phi(w)\circ\Phi(s)^{-1}=\overline\Phi(f)\).
In particular, a W-morphism \(b=(r,f):x\to y\) gives the cone
\((x,\id_x,b)\). Applying \eqref{eq:roofindependence} to this
cone yields
\[
 \overline\Phi(f)=\Phi(b).
\]

This includes identity morphisms and proves agreement with \(\Phi\)
on every existing square. It also shows that all W-witnesses for a
fixed value function have the same image under \(\Phi\).
For composition, let \(f:x\to y\) and \(g:y\to v\) be
AA-morphisms. Proposition~\ref{prop:ore}, applied to \(w_f\) and
\(s_g\), supplies \(t:u\to z_f\) in \(\Sts\) and a
W-morphism \(h:u\to z_g\) such that
\(w_f\circ t=s_g\circ h\). Thus
\(\Phi(s_g)^{-1}\circ\Phi(w_f)
=\Phi(h)\circ\Phi(t)^{-1}\), and
\begin{align*}
 \overline\Phi(g)\circ\overline\Phi(f)
 &=\Phi(w_g)\circ\Phi(s_g)^{-1}
       \circ\Phi(w_f)\circ\Phi(s_f)^{-1}\\
 &=\Phi(w_g\circ h)\circ\Phi(s_f\circ t)^{-1}\\
 &=\overline\Phi(g\circ f).
\end{align*}

The last equality follows from \eqref{eq:roofindependence}, since
\((u,s_f\circ t,w_g\circ h)\) represents \(g\circ f\).
Consequently \(\overline\Phi\) is a functor and
\(\overline\Phi\circ Q=\Phi\).
In \(\NetAA\), the graph cone satisfies
\(f=Q(w_f)\circ Q(s_f)^{-1}\). Any factorization of \(\Phi\)
must therefore satisfy \eqref{eq:extension}. Its object values are
also forced because \(Q\) is the identity on objects. This proves
uniqueness of the factorization.

It remains to treat natural transformations. Let
\(\alpha:\Phi\Rightarrow\Psi\) be a natural transformation between
two functors inverting \(\Sts\). Define \(\overline\alpha_x=\alpha_x\).
Naturality for \(s_f\) and \(w_f\) gives
\begin{align*}
 \overline\Psi(f)\circ\alpha_x
 &=\Psi(w_f)\circ\Psi(s_f)^{-1}\circ\alpha_x\\
 &=\Psi(w_f)\circ\alpha_{z_f}\circ\Phi(s_f)^{-1}\\
 &=\alpha_y\circ\Phi(w_f)\circ\Phi(s_f)^{-1}
  =\alpha_y\circ\overline\Phi(f).
\end{align*}

Hence these same components define a natural transformation
\(\overline\alpha:\overline\Phi\Rightarrow\overline\Psi\).
The components force uniqueness, and their unchanged values ensure
preservation of identities and vertical composition of natural
transformations. Extension followed by restriction is the identity.
Conversely, extending the restriction of a functor on \(\NetAA\)
recovers that functor by the graph formula. Thus \(Q^*\) is an
isomorphism of categories for \(S=\Sts\).

Finally, a functor inverting \(\Seq\) or \(\Ssat\) also
inverts \(\Sts\) and has the factorization just proved. Every
functor factoring through \(Q\) inverts both larger classes.
The same restriction and extension operations therefore give the
isomorphism for either of them.
\end{proof}

The theorem identifies \(\Ssat\) as the saturation of \(\Sts\)
and of \(\Seq\). Indeed, each morphism of \(\Ssat\) is sent to
an isomorphism by every functor inverting either \(\Sts\) or
\(\Seq\), because every such functor factors through \(Q\). A morphism outside
\(\Ssat\) is not forced to become invertible, as \(Q\) itself
shows. Also, parallel morphisms in \(\Wpar\) have equal images in
this localization exactly when their value functions agree.

A Kelley index map need not be monotone. The next observation changes
the order on its domain so that it becomes a Willard index map, while
preserving the eventuality filter of the source net.

\begin{lemma}\label{lem:korder}
Let \(r:\mathbb D\to\mathbb E\) be eventually cofinal. Define
a preorder on the underlying set of \(\mathbb D\) by
\[
 d\le_r d'\quad\Leftrightarrow\quad
 d\le d'\text{ in }\mathbb D
 \text{ and }r(d)\le r(d')\text{ in }\mathbb E.
\]

Then \(\mathbb D_r=(\mathbb D,\le_r)\) is directed and has
filter \(\TT_{\mathbb D}\). The identity map
\(a:\mathbb D_r\to\mathbb D\) is monotone and tail-surjective,
and \(r:\mathbb D_r\to\mathbb E\) is monotone cofinal.
\end{lemma}
\begin{proof}
For \(d,d'\in\mathbb D\), choose a common upper bound of
\(r(d),r(d')\) in \(\mathbb E\). Eventual cofinality gives
a tail set of \(\mathbb D\) whose images lie above that bound.
An element of this tail set above \(d,d'\) is a common upper bound
in \(\mathbb D_r\), proving directedness.

For a fixed \(d_0\), choose \(d_*\) such that
\(r(d_0)\le r(d)\) whenever \(d\ge d_*\), and choose
\(d_1\ge d_0,d_*\). Then
\[
 \mathbb D_{\ge d_1}
 \subseteq(\mathbb D_r)_{\ge_r d_0}
 \subseteq\mathbb D_{\ge d_0}.
\]

These inclusions prove equality of the two filters generated by tail sets and
tail-surjectivity of \(a\). The map \(a\) is monotone by
definition of \(\le_r\). The same definition makes \(r\)
monotone, and its cofinal range is unchanged.
\end{proof}

We now extend the localization to Kelley morphisms and to the categories
in which a morphism records only the value function, while requiring an
index map to exist.
In \(\Kpar\), we keep the notation \(\Sts\), \(\Seq\) and \(\Ssat\)
for the Willard morphisms in \(\Sts\), the exact vertical
K-morphisms and \(Q_{\Kpar}^{-1}(\Iso(\NetAA))\), respectively.
In \(\Wex\) and \(\Kex\), a morphism belongs to one of
these classes when its value function admits a witness from the
corresponding class of \(\Wpar\) or \(\Kpar\).

In \(\Wex\), the index map \(n\mapsto2n\) of Example~\ref{ex:classes}
no longer separates the first two classes, since its value function
also admits the identity witness. Nevertheless, the two classes still
differ. Consider one-point constant nets whose index sets are \(\NN\)
and \(\NN^2\), respectively. The diagonal is a W-witness for an exact
morphism from the first net to the second. Every monotone map
\(\NN\to\NN^2\) has chain image, whereas every product tail set
contains incomparable points \((p+1,q)\) and \((p,q+1)\).
No such witness can be tail-surjective.

\begin{theorem}\label{thm:fourlocal}
Let \(\CC\) be any of
\(\Wpar,\Wex,\Kpar,\Kex\). The functor
\(Q_{\CC}:\CC\to\NetAA\) is a localization at each of the classes
\(\Sts\), \(\Seq\) and \(\Ssat\) of \(\CC\). More precisely, for
every category \(\mathsf A\), the functor in \eqref{eq:universal} is
an isomorphism. In particular,
\(\CC[\Seq^{-1}]\simeq\NetAA\simeq\Filp\).
\end{theorem}
\begin{proof}
For \(\Kpar\), suppose \(\Phi\) inverts \(\Sts\). Restrict
it to \(\Wpar\) and apply Theorem~\ref{thm:wlocal} to obtain
\(\overline\Phi:\NetAA\to\mathsf A\). Given a K-morphism
\(b=(r,f):x\to y\), regard its source as the same function
\(x_r:\mathbb D_r\to X\) with the order in
Lemma~\ref{lem:korder}. That lemma gives the W-morphism
\(a=(\id_{\mathbb D},\id_X):x_r\to x\), which belongs to
\(\Sts\), and the W-morphism \(t=(r,f):x_r\to y\), with
\(b\circ a=t\). Therefore
\[
 \Phi(b)=\Phi(t)\circ\Phi(a)^{-1}
        =\overline\Phi(Q_{\Kpar}(b)).
\]

The last equality uses the already established factorization on
W-morphisms and the identity
\(Q_{\Kpar}(b)=Q(t)\circ Q(a)^{-1}\). Thus
\(\overline\Phi\) factors \(\Phi\) on all K-morphisms, and
uniqueness follows by restriction to \(\Wpar\).
A natural transformation between two such functors extends uniquely
from their Willard restrictions by Theorem~\ref{thm:wlocal}.
The resulting transformation is natural for every AA-morphism, and its
restriction to the Kelley category has the original components.
This proves the isomorphism of functor categories for \(\Kpar\).
The larger classes map to isomorphisms by Lemma~\ref{lem:iso}, so
the enlargement argument of Theorem~\ref{thm:wlocal} applies.

For the remaining cases, let \(P=P_W\) or \(P=P_K\), and
write \(\CC\) and \(\CC_\exists\) for its source
and target, respectively. If \(\Phi:\CC_\exists\to\mathsf A\) inverts
the indicated denominators, then \(\Phi\circ P\) inverts their
corresponding class in \(\CC\). The result just proved
gives a unique \(\overline\Phi:\NetAA\to\mathsf A\) with
\[
 \overline\Phi\circ Q_{\CC}=\Phi\circ P.
\]

Every morphism \(f\) of \(\CC_\exists\) has a lift \(b\) in
\(\CC\), because \(P\) is full. Evaluation at this
lift gives
\(\overline\Phi(Q_{\CC_\exists}(f))=\Phi(f)\).
The equality on objects holds because \(P\) is the identity on
objects, so the factorization is strict throughout.

For a natural transformation \(\alpha:\Phi\Rightarrow\Psi\),
precompose with \(P\) and extend by the result for
\(\CC\). Its extended components remain \(\alpha_x\).
Naturality for each morphism of \(\CC_\exists\) follows by choosing
a lift, and uniqueness follows by precomposing with \(P\).
Restriction and extension are therefore inverse on both functors
and natural transformations.
\end{proof}

\begin{corollary}\label{cor:filterslocal}
For each \(\CC\in\{\Wpar,\Wex,\Kpar,\Kex\}\), the composite
\[
 \CC\overset{Q_{\CC}}{\rightarrow}\NetAA\overset{T}{\rightarrow}\Filp
\]
is a localization at each corresponding denominator class specified
in Theorem~\ref{thm:fourlocal}.
\end{corollary}
\begin{proof}
Precomposition with \(Q_{\CC}\) is an equivalence by
Theorem~\ref{thm:fourlocal}, and precomposition with \(T\)
is an equivalence because \(T\) is an equivalence. Their composite is the required
equivalence for \(T\circ Q_{\CC}\).
\end{proof}

\begin{remark}\label{rem:germs}
The morphisms of \(\Filp\) in Corollary~\ref{cor:filterslocal}
are functions, and equality of morphisms is ordinary equality of
functions. A different comparison arises when functions are allowed
to have domains restricted to members of a filter. A \emph{germ}
relative to a filter \(\FF\) is represented by such a function,
defined on some \(A\in\FF\). Two representatives determine
the same germ when they agree on a member of \(\FF\) contained in
both domains.

Blass~\cite[Section 1 and Theorem 15]{Blass}, allowing also improper
filters, obtains his category of filters and germs by localizing the
category of filters and functions at inclusions of members
of the filters. For ultrafilters,
Garner~\cite[Section 2.2, Definition 4 and Proposition 5]{Garner}
describes this passage by inverting the continuous injections.
These localizations make further identifications beyond those in
Corollary~\ref{cor:filterslocal}.

For example, consider the constant zero net in \(X=\{0,1\}\).
Both \(\id_X\) and the constant zero function give W-endomorphisms
with the identity index map. They agree on \(\{0\}\), a member of
the eventuality filter, and hence determine the same germ. They
remain distinct in \(\NetAA\) and \(\Filp\), since their values
at \(1\) differ.

These localizations also make further morphisms invertible. Let
\(i:\{0\}\to X\) be the inclusion. With the identity index map, it
gives a W-morphism from the constant zero net in \(\{0\}\) to the
constant zero net in \(X\). Since \(i\) is not bijective, its image
is not an isomorphism in \(\NetAA\) by Lemma~\ref{lem:iso}, nor in
\(\Filp\). It is, however, the inclusion of the member \(\{0\}\) of
the eventuality filter, and it becomes invertible after the passage
to germs, with inverse represented by the constant function
\(X\to\{0\}\).
\end{remark}

The existence of the localization does not by itself imply that its
denominators admit a calculus of fractions. Recall that such a
calculus describes localized morphisms by pairs of morphisms in the
original category, together with rules for composing and comparing
these pairs. We use the convention of a \emph{right multiplicative
system} in~\cite[Definition 4.27.1]{Stacks}. The pairs have the form
\[
\begin{tikzcd}[column sep=2em,row sep=1.5em]
& u \arrow[dl,"s"'] \arrow[dr,"h"] & \\
x & & y
\end{tikzcd}
\]
with \(s\in S\), and represent the composite of the image of
\(h\) with the inverse of the image of \(s\) in the localization.
The required conditions on \(S\) are:
\begin{enumerate}[label=(F\arabic*)]
\item \(S\) contains identities and is closed under composition.
\item For \(f:x\to y\) and \(s:z\to y\) with \(s\in S\),
there are \(t:u\to x\) in \(S\) and \(h:u\to z\) such that
\(f\circ t=s\circ h\).
\item If \(s\circ f=s\circ g\) with \(s\in S\), there is
\(t\in S\), with target the common source of \(f,g\), such that
\(f\circ t=g\circ t\).
\end{enumerate}

Proposition~\ref{prop:ore} establishes (F2) for \(\Sts\) in
both \(\Wpar\) and \(\Kpar\), and (F1) was checked above.
Condition (F3), however, fails. Enlarging \(\Sts\) to \(\Seq\)
or \(\Ssat\) also causes (F2) to fail, as the following maps
show.
Take the constant net \(c:\NN\to\{*\}\). The vertical
endomorphisms with index maps \(\rho_0(n)=2n\) and \(\rho_1(n)=2n+1\)
belong to \(\Seq\). A square as in (F2) would require
\(\rho_0\circ a=\rho_1\circ b\), which is impossible on a nonempty
domain because the two ranges are disjoint. Thus (F2) fails for
\(\Seq\) and \(\Ssat\).

The vertical endomorphism with index map
\(\lambda(n)=\lfloor n/2\rfloor\) belongs to \(\Sts\), and
\(\lambda\circ\rho_0=\lambda\circ\rho_1=\id_{\NN}\). Nevertheless,
\(\rho_0\circ v=\rho_1\circ v\) is impossible on any nonempty domain.
Hence (F3) fails for \(\Sts\) and for the two larger classes.
All these maps are W-morphisms, so the same obstructions apply to the
corresponding classes in \(\Kpar\).

\begin{proposition}\label{prop:fractions}
In either \(\Wex\) or \(\Kex\), each of the three denominator
classes specified in Theorem~\ref{thm:fourlocal} is a right
multiplicative system. For the smallest class in \(\Kex\), the
condition is admission of a tail-surjective Willard witness.
\end{proposition}
\begin{proof}
Let \(\CC\) be \(\Wex\) or \(\Kex\), and let \(S\) be one
of the indicated classes. For (F1), identities and composites admit
the required witnesses whenever the class is defined by such an
existence condition. Exactness and the preimage of isomorphisms are
also closed under composition. All three classes contain the image
of \(\Sts\) and are sent to isomorphisms by \(Q_{\CC}\).

Given \(f:x\to y\) and \(s:z\to y\) with \(s\in S\),
the composite \(Q_{\CC}(s)^{-1}\circ Q_{\CC}(f):x\to z\)
is an AA-morphism. Theorem~\ref{thm:graph} represents its value
function by a morphism \(t:u\to x\) admitting a witness in
\(\Sts\) and a W-morphism \(h:u\to z\). The value functions
satisfy \(f\circ t=s\circ h\), which is equality of morphisms
in \(\CC\). This proves (F2).

For (F3), Lemma~\ref{lem:iso} shows that the value function of
\(s\) is bijective. Therefore \(s\circ f=s\circ g\) implies
\(f=g\) as functions and hence as morphisms of \(\CC\).
Taking \(t\) to be the identity proves the condition.
\end{proof}

In \(\Wpar\) and \(\Kpar\), the index witnesses obstruct (F3)
for all three classes and (F2) for the two larger ones. Once only the
existence of witnesses is retained, all three classes become right
multiplicative systems, while Theorem~\ref{thm:fourlocal} shows that
the localization is the same. Thus, in \(\Wex\) and \(\Kex\),
AA-morphisms are described by pairs \((s,h)\) as above with the
standard rules of the calculus of fractions~\cite[Section 4.27]{Stacks},
whereas in \(\Wpar\) and \(\Kpar\) the representations of
Section~\ref{sec:graphs} have to be used instead.

\section{Final remarks}\label{sec:final}
The representations of Section~\ref{sec:graphs} can be read as
fractions. Let \(\CC\) be
\(\Wpar\) or \(\Kpar\), and let \((u,s,w)\) represent an AA-morphism
\(f:x\to y\). The denominator \(s\in\Sts\) changes the presentation
of \(x\) without changing its eventuality filter, while the numerator
\(w\) has value function \(f\). In \(\NetAA\),
\(f=Q_{\CC}(w)\circ Q_{\CC}(s)^{-1}\). The terminal cone
\((z_f,s_f,w_f)\) of Theorem~\ref{thm:graph} is analogous to a reduced
fraction: every cone \((u,s,w)\) representing \(f\) has a unique
morphism \(m:u\to z_f\) with \(s=s_f\circ m\) and
\(w=w_f\circ m\). Similarly, every fraction of integers with positive
denominator is obtained from its reduced form by multiplying numerator
and denominator by the same positive integer. The
completion in Proposition~\ref{prop:ore} plays a role analogous to
finding a common denominator when these representations are composed.
As one would expect, the analogy has limitations: the order of the
factors matters, and the morphism \(m\) need not belong to \(\Sts\),
although \(Q_{\CC}(m)\) is invertible. Moreover, \(\Sts\)
fails condition (F3) in \(\Wpar\) and \(\Kpar\), so it does not give a
calculus of right fractions in these categories.

We now return to the comparison of nets and filters on one fixed
set \(X\). For any of the four categories \(\CC\) in
Theorem~\ref{thm:fourlocal}, let \(\CC_X\) have the nets in \(X\) as
objects and, as morphisms, those of \(\CC\) whose value function is \(\id_X\).
The value-set functor \(\CC\to\mathbf{Set}\) sends a net to the
set in which it takes values and a morphism to its value function.
So \(\CC_X\) is simply its fiber over \(X\).\footnote{For
\(\CC=\Kpar\), this fiber is the comma category \((F\downarrow X)\)
of \v{S}lapal~\cite[Section 2]{Slapal2016} recalled in
Section~\ref{sec:categories}.}
Let \(\mathsf P_X\) be the category of proper filters on \(X\), with
one morphism \(\FF\to\GG\) precisely when \(\GG\subseteq\FF\), and let
\(S_X\) consist of the morphisms in \(\CC_X\) between nets with the
same eventuality filter.

\begin{samepage}
\begin{corollary}\label{cor:fixedset}
The functor \(E_X:\CC_X\to\mathsf P_X\), given on objects by
\(E_X(x)=\Ev(x)\), is a localization at \(S_X\). In particular,
\[
 \CC_X[S_X^{-1}]\simeq\mathsf P_X.
\]

The same conclusion holds for the smaller denominator class consisting
of morphisms admitting a tail-surjective Willard witness, with that
witness required to be the specified index map in a category of squares.
\end{corollary}
\end{samepage}
\begin{proof}
Let \(\mathsf A_X\) have the nets in \(X\) as objects and a unique
morphism \(x\to y\) whenever \(\Ev(y)\subseteq\Ev(x)\). All the
constructions used to prove Theorem~\ref{thm:fourlocal} preserve the
value set when the value functions are identities: this holds for
the representation and comparison in Theorem~\ref{thm:graph}, the
completion in Proposition~\ref{prop:ore}, and the order refinement
in Lemma~\ref{lem:korder}. The universal-property proof therefore
restricts to \(\CC_X\to\mathsf A_X\), including its assertion about
natural transformations. It gives localization at the smaller class
and, by the same enlargement argument, at \(S_X\).

The functor \(\mathsf A_X\to\mathsf P_X\) taking eventuality filters
is fully faithful by the definition of its hom-sets, and essentially
surjective by Proposition~\ref{prop:base}. Composing with this
equivalence proves the assertion.
\end{proof}

Consequently the proper filters on \(X\), ordered by refinement,
form a localization of nets in that set as well, and the reverse
inclusion in \(\mathsf P_X\) records the direction of the subnet
relation: a subnet has a finer eventuality filter. Although the nets
in \(X\) may form a proper class, the localization is equivalent to
this small poset.\footnote{This description is also compatible with
changes of value set. A function \(f:X\to Y\) sends \(x\) to
\(f\circ x\), keeps the index map of each vertical square, and sends
\(\FF\) to \(f(\FF)\). These assignments preserve the denominator
classes and satisfy \(\Ev(f\circ x)=f(\Ev(x))\).}

The representation theorem also lets us lift a subnet through a
function. Let \(f:A\to B\) be a function. For sequences
and subsequences, the construction is immediate: if
\(\varphi:\mathbb N\to A\) and
\(\psi=(f\circ\varphi)\circ\rho\), where
\(\rho:\mathbb N\to\mathbb N\) is strictly increasing, then
\(\eta=\varphi\circ\rho\) is a subsequence of \(\varphi\)
and \(f\circ\eta=\psi\). The maps form the commutative diagram
\[
\begin{tikzcd}[row sep=large,column sep=huge]
\mathbb N\arrow[r,"\varphi"]&A\arrow[d,"f"]\\
\mathbb N\arrow[u,"\rho"]\arrow[ur,"\eta",description]
  \arrow[r,"\psi"']&B.
\end{tikzcd}
\]

For an AA-subnet of an image net, a map such as \(\rho\) need
not be available. The next corollary supplies a lift after passing
to a W-subnet of the given AA-subnet, while preserving its
eventuality filter.

\begin{corollary}\label{cor:subnetlifting}
Let \(f:A\to B\) be a function, \(\varphi:\mathbb D\to A\) be
a net, and \(\psi:\mathbb E\to B\) be an AA-subnet of
\(f\circ\varphi\). There is a W-subnet \(\eta\) of \(\varphi\)
such that \(f\circ\eta\) is a W-subnet of \(\psi\) and
\[
 \Ev(f\circ\eta)=\Ev(\psi).
\]
\end{corollary}
\begin{proof}
Apply Theorem~\ref{thm:graph} to the AA-morphism
\(\id_B:\psi\to f\circ\varphi\). It gives a net
\(\theta:\mathbb C\to B\) and monotone cofinal maps
\(q:\mathbb C\to\mathbb E\), \(r:\mathbb C\to\mathbb D\) such
that
\[
 \theta=\psi\circ q=f\circ\varphi\circ r.
\]

Moreover, \(\Ev(\theta)=\Ev(\psi)\). Set \(\eta=\varphi\circ r\). The maps in the construction form
the commutative diagram
\[
 \begin{tikzcd}[row sep=large,column sep=huge]
  \mathbb D\arrow[r,"\varphi"]&A\arrow[d,"f"]\\
  \mathbb C\arrow[u,"r"]\arrow[ur,"\eta",description]
    \arrow[r,"\theta"']\arrow[d,"q"']&B.\\
  \mathbb E\arrow[ur,"\psi"']&
 \end{tikzcd}
\]

Both index maps are monotone cofinal, and \(q\) is tail-surjective.
Since \(f\circ\eta=\theta\), the required assertions follow.
\end{proof}

This lifting result also bears on compactness, through ultranets.
Recall that an \emph{ultrafilter} is a proper filter
maximal under inclusion, and that a net \(u\) in \(X\) is an
\emph{ultranet} when \(\Ev(u)\) is an ultrafilter.\footnote{The
authors consider \emph{infranet} a more suitable term: the subnet
characterization below expresses minimality, up to tail-equivalence,
for the preorder defined by the AA-subnet relation.}
By Proposition~\ref{prop:base}, every proper filter containing
\(\Ev(u)\) is the eventuality filter of an AA-subnet of \(u\), so
\(u\) is an ultranet exactly when it is tail-equivalent to each of its
AA-subnets. Corollary~\ref{cor:subnetlifting} then shows that functions
preserve ultranets. If \(\varphi\) is an ultranet and \(\psi\) is an
AA-subnet of \(f\circ\varphi\), the W-subnet \(\eta\) given by the
corollary satisfies \(\eta\sim\varphi\), since every subnet of an
ultranet is tail-equivalent to it. Hence
\(\Ev(\psi)=\Ev(f\circ\eta)=f(\Ev(\varphi))=\Ev(f\circ\varphi)\), so
\(f\circ\varphi\) is tail-equivalent to each of its AA-subnets.
This is the preservation step in the ultranet proof of Tychonoff's
theorem, which uses that a space is compact if and only if every
ultranet in it converges~\cite[Section~7.5.6]{MuegerTopology}. There
the step is usually obtained from inverse images of subsets. Here it
follows from the representation used to prove the localization.

One might compare the passage from nets to filters with the
identification of measurable functions that agree almost everywhere.
We expect our work to show how this comparison is misleading in at
least three respects. First, the indices of
a net parametrize a base of its eventuality filter, and passing to the
filter forgets this parametrization, much as describing the circle as
a subset of the plane forgets its parametrization by angles. Regarding
nets as dispensable for this reason would be like regarding convergent
series of real numbers as dispensable because one can work directly
with their sums. Second, nets come with evident morphisms, the
commutative squares formed by index maps and value functions, and the
localization describes how these morphisms give rise to those of
\(\Filp\). Third, the passage involves more than the identification
of tail-equivalent nets. For nets in a fixed set, the AA-subnet
relation is a preorder, with tail-equivalence as its symmetric part.
The map \(\Ev\) reverses this preorder with respect to inclusion:
\(x\) is an AA-subnet of \(y\) exactly when
\(\Ev(y)\subseteq\Ev(x)\). The analogue
of agreement almost everywhere is rather the passage to germs, which
identifies functions agreeing on a member of a filter. As
Remark~\ref{rem:germs} illustrates, this is a further step that
\(\Filp\) does not take.

The equivalence therefore does not make nets dispensable. That would
require every relevant use of a net to depend only on its eventuality
filter, and the equivalence does not say this. Constructions may use
indices explicitly, as the lifting argument above does. What the
localization says is which changes of presentation a construction
must ignore in order to descend to filters.

\bigskip
\noindent\begin{minipage}{\textwidth}
\small
\textsc{Germ\'an I. Gomero Ferrer}\\
Departamento de Ci\^encias Exatas, Universidade Estadual de Santa Cruz\\
Ilh\'eus, BA, 45662-900, Brazil\\
\textit{Email address:} \href{mailto:gigferrer@uesc.br}{\nolinkurl{gigferrer@uesc.br}}
\end{minipage}

\medskip
\noindent\begin{minipage}{\textwidth}
\small
\textsc{Renan M. Mezabarba}\\
Departamento de Ci\^encias Exatas, Universidade Estadual de Santa Cruz\\
Ilh\'eus, BA, 45662-900, Brazil\\
\textit{Email address:} \href{mailto:rmmezabarba@uesc.br}{\nolinkurl{rmmezabarba@uesc.br}}
\end{minipage}

\end{document}